\documentclass[11pt, reqno]{amsart}

\usepackage{color}

 \usepackage{amsmath}
 \usepackage{amssymb}
\usepackage{bm}

\newcommand{\mysection}[1]{\section{#1}
      \setcounter{equation}{0}}

\newtheorem{theorem}{Theorem}[section]
\newtheorem{lemma}[theorem]{Lemma}

\newtheorem{corollary}[theorem]{Corollary} 

\theoremstyle{definition}

\theoremstyle{remark}
\newtheorem{remark}[theorem]{Remark}

\newcommand{\loc}{\text{\rm loc}}

 \makeatletter 
 \def\dashint{%
 \operatorname%
 {\,\,\text{\bf--}\kern-.98em\DOTSI\intop\ilimits@\!\!}}
\def\dashnorm{\,\,\text{\bf--}\kern-.5em\|}
\def\ninf{\qopname\relax\@empty{inf\phantom{p}\!\!\!}}

 \makeatother

\newcommand\bbeta{\text{\raise-.2ex\hbox{$\bm{\beta}$}}}

\newcommand\bR{\mathbb{R}}

\newcommand\cC{\mathcal{C}}

\begin{document}

\title[Parabolic Adams's,   Chiarenza-Frasca  theorems]
{On parabolic Adams's, the  Chiarenza-Frasca
theorems,
and some other results related to parabolic Morrey spaces}

\author{N.V. Krylov}
 
\email{nkrylov@umn.edu}
\address{127 Vincent Hall, University of Minnesota,
 Minneapolis, MN, 55455}



\begin{abstract}
We present several results of embedding type
for parabolic Morrey and $L_{p}$ spaces with or without
mixed norms. Some other interpolation
results for parabolic Morrey spaces are also given. The main object
of investigation is the term $b^{i}D_{i}u$
and the ways to estimate it in various Morrey 
and $L_{p}$ spaces in order to be able to treat it as a perturbation term in the parabolic
equations.

\end{abstract}

\maketitle

\mysection{Introduction}

Let $\bR^{d}$, $d\geq 1$, be a Euclidean space of
points $x=(x^{1},...,x^{d})$.
In 1975 D.~Adams \cite{Ad_75} among many other things proved that, if $d\geq2 $ and we are given $u\in C^{\infty}_{0}
=C^{\infty}_{0}(\bR^{d})$ with its gradient
$Du=(D_{1}u,...,D_{d}u)$, $D_{i}=\partial/\partial x^{i}$, satisfying
\begin{equation}
                             \label{1.30.2}
\int_{|x-y|<\rho}|Du(y)|^{q}\,dy\leq  \rho^{d-\beta q},
\end{equation}
with $q>1$, $1<\beta\leq d/q$, and any $\rho\in(0,\infty)$  and
$x\in \bR^{d}$, then for all $\rho\in(0,\infty)$ and
$x\in \bR^{d}$ we have
\begin{equation}
                             \label{1.30.1}
\int_{|x-y|<\rho}| u(y)|^{r}\,dy\leq N\rho^{d-(\beta-1)r}
\end{equation}
with a   constant $N$ independent of $u$ and $r$ satisfying $(\beta-1)r=\beta q$.

This fact played a crucial role in \cite{Kr_303} where the author investigated the
solvability of elliptic equations
\begin{equation}
                                \label{1.29.1}
a^{ij}D_{ij}u+b^{i}D_{i}u+u=f \quad (D_{ij}=D_{i}D_{j})
\end{equation}
with $b\not \in L_{d,\loc}$ but rather satisfying for a sufficiently small $\hat b$,
  all sufficiently small $\rho$ and all balls
$B$ of radius $\rho$
$$
\int_{B}|b|^{d_{0}}\,dx
\leq \hat b \rho^{d-d_{0}}
$$
with certain $d_{0}\in(d/2,d)$.

Our goal in this paper is to prepare necessary
tools for developing a similar theory for parabolic equations. In Section
\ref{section 10.8.1} we prove an analog
of Adams's intermediate estimate, which
is the main starting point. Section
\ref{section 10.8.2} contains the parabolic
analog of the embedding theorem mentioned in the beginning of the article. It also contains ``local'' 
interpolation  inequalities in Morrey
spaces allowing one to deal with Morrey's norms of expressions
like $b^{i}D_{i}u$ in domains when $b$ is bounded.
Section \ref{section 1.25.1} is devoted
to the parabolic analog of a Chiarenza-Frasca
theorem allowing to estimate the $L_{p}$-norm
rather than Morrey's norm
of $b^{i}D_{i}u$. In Section
\ref{section 10.8.3} we treat parabolic Morrey spaces with mixed norms. The main object
of investigation is the term $b^{i}D_{i}u$
and the ways to estimate it in various Morrey 
and $L_{p}$ spaces in order to be able to treat it as a perturbation term in the parabolic
analog of \eqref{1.29.1}.

We finish the introduction with some notation
and a remark.
Define $B_{\rho}(x)=\{y \in\bR^{d} :|x-y|<\rho\}$,
 $\bR^{d+1}:=
\{z=(t,x):t\in\bR,x\in\bR^{d}\}$,
$$
C_{\rho}(t,x)=\big\{(s,y)\in \times\bR^{d +1 } :|x-y|<\rho,
t\leq s< t+\rho^{2}\big\},\quad C_{\rho}=C_{\rho}(0)
$$
and let $\cC_{\rho}$ be the collection of
$C_{\rho}(z)$, $z\in\bR^{d+1}$, $\cC=
\{\cC_{\rho},\rho>0\}$. 
For measurable $\Gamma\subset \bR^{d+1}$
set $|\Gamma|$ to be its Lebesgue measure
and when it makes sense set
$$
f_{\Gamma}=\dashint_{\Gamma}f\,dz=\frac{1}{|\Gamma|}
\int_{\Gamma}f\,dz.
$$
Similar notation is used for $f=f(x)$.

\begin{remark}
Formally, Adams proved \eqref{1.30.1}  
assuming that $d\geq2$. However, it is also true if $d=1$. To show this it suffices
to take $u$ depending only on one coordinate.
The reader may wonder how the restriction
$\beta\leq d/q$ will become $\beta\leq 1/q$.
The point is that if $d=1$ and $\beta> 1/q$,
we have $d-\beta q<0$ and  condition
\eqref{1.30.2} becomes only possible if $u=0$.

\end{remark}
\mysection{Preliminary estimates}

                       \label{section 10.8.1}
An important quantity characterizing $L_{p}
=L_{p}(\bR^{d+1})$ is what we call the index
which is the exponent of $\rho$ in the expression
$$
\| I_{C_{\rho}}\|_{L_{p}}\quad
\text{\rm that is}\quad \frac{d+2}{p} .
$$

For domains $Q\subset\bR^{d +1 }$, $p\in[1,\infty)$, and $\beta\in(0,(d+2)/p]$, introduce
Morrey's space $E_{p,\beta}(Q)$
as the set of $g $ such that
\begin{equation}
                             \label{8.11.2}
\|g\|_{E_{p,\beta}(Q)}:=
\sup_{\rho<\infty,(t,x)\in Q}\rho^{\beta}
\dashnorm g I_{Q}\|_{L_{p}(C_{ \rho}(t,x) )} <\infty ,
\end{equation}
where
$$
\dashnorm g\|_{L_{p}(\Gamma)}=\Big(
\dashint_{\Gamma}|g|^{p}\,dz\Big)^{1/p}.
$$
We abbreviate $E_{p,\beta}=E_{p,\beta}(\bR^{d+1})$.
Observe that if $Q=C_{R}$
one can restrict $\rho$ in \eqref{8.11.2}
to $\rho\leq R$ since $\beta\leq (d+2)/p$. 
Also in that case one can allow $(t,x)$ to
be arbitrary, because, if $|x|\geq R$,
then $B_{R}\cap B_{\rho}(x)\subset
 B_{R}\cap B_{\rho}(Rx/|x|)$.
It is also useful to observe that,
in case $Q=C_{R}$, one gets an equivalent norm
by adding to the restrictions $\rho<\infty,
(t,x)\in C_{R}$, the requirement that the geometric
center of $C_{\rho}(t,x)$ be in $C_{R}$. This follows from the fact that the $L_{p}(C_{\rho}(t,x))$-norm of $gI_{C_{R}}$ will only increase
if we pull $C_{\rho}(t,x) $ down
the $t$ axis to $\{t=0\}$ (if $\rho^{2}>2R^{2} $) or to the moment that the shifted $C_{\rho}(t,x) $ has its geometric center inside $C_{R}$.

There are many different notations for
the norms in Morrey spaces. The convenience
of the above notation is well illustrated
by Theorem \ref{theorem 9.29.1}
and Corollary \ref{corollary 10.8.1}.
 
We will often, always tacitly, use the following formulas
in which $u(t,x)=v(t/R^{2},x/R)$:
$$
\dashnorm u\|_{L_{p}(C_{R})}
=\dashnorm v\|_{L_{p}(C_{1})},\quad
\|u\|_{E_{p,\beta}(Q)}=R^{\beta}
\|v\|_{E_{p,\beta}(Q_{R})}, 
$$
where $Q_{R}=\{(t,x): (  R^{2}t,Rx )\in Q\}$,
$$
\|Du\|_{E_{p,\beta}(C_{R})}=R^{\beta-1}
\|Dv\|_{E_{p,\beta}(C_{1})},\quad
\|D^{2}u\|_{E_{p,\beta}(C_{R})}=R^{\beta-2}
\|v\|_{E_{p,\beta}(C_{1})}.
$$

For $s,r>0,\alpha>0$, and appropriate $f(t,x)$'s
on $\bR^{d+1}$ define
$$
p_{\alpha}(s,r)=\frac{1}{s^{(d+2-\alpha)/2}}e^{-r^{2}/s}I_{s>0},
$$
$$
P_{\alpha}f(t,x)=\int_{\bR^{d+1} }p_{\alpha}(s,|y|)f(t+s,x+y)\,dyds.
$$
Observe that, if $f$ is independent of $t$,
then
$$
P_{\alpha}f(t,x)=
P_{\alpha}f(x)=N(\alpha)\int_{\bR^{d}}
\frac{1}{|y|^{d-\alpha}}f(x+y)\,dy= NI_{\alpha}f(x),
$$
where $I_{\alpha}$ is the Riesz potential.
Therefore, one can get the Adams 
estimate found in the proof of Proposition 3.1
of \cite{Ad_75} from \eqref{9.28.1} below.
In our investigation the most important values
of $\alpha$ are 1 and 2.
Set
$$
M_{\beta}f(t,x)=\sup_{\rho>0}\rho^{\beta}
\dashint_{C_{\rho}(t,x)}|f(z)|\,dz,\quad
0\leq\beta\leq d+2,
$$
$$
M f =M_{0}f .
$$
The following lemma is obtained by integrating by parts.
\begin{lemma}
                      \label{lemma 1.16.1}
Let $\beta>0$ be a finite number, $f(t)\geq0$ be a   function on $[0,\infty)$ such that
$$
t^{-\beta}\int_{0}^{t}f(s)\,ds\to 0 
$$
as $t\to \infty$. Then, for any $  S\geq 0$,
$$
\int_{S}^{\infty}t^{-\beta}f(t)\,dt\leq
\beta \int_{S}^{\infty}t^{-\beta-1}\Big(\int_{S}^{t}f(s)\,ds\Big)\,dt.
$$

\end{lemma}

\begin{lemma}
                           \label{lemma 9.28.1}
For any $\alpha\in(0,\beta ),\beta\in(0,d+2]$ there exist   constants $N$
($<\infty$)  such that for any $f\geq0$
and $\rho\in(0,\infty)$ we have
\begin{equation}
                             \label{1.17.2}
P_{\alpha}(I_{C _{ \rho}}f)(0)
\leq N\rho^{\alpha }M 
f(0) ,\quad
P_{\alpha}(I_{C^{c}_{ \rho}}f)(0) 
\leq N\rho^{\alpha-\beta}M_{\beta}
f (0),
\end{equation}
\begin{equation}
                                \label{9.28.1}
 P_{\alpha}f \leq    
N(M_{\beta}f )^{\alpha/\beta}
(  M f)^{1-\alpha/\beta}.
\end{equation}
In particular (by H\"older's inequality), for any $p\in[1,\infty]$,
$q\in(1,\infty]$, and measurable $\Gamma$
\begin{equation}
                                \label{9.28.2}
 \|P_{\alpha}f\|_{L_{r}(\Gamma)} \leq    
N\|M_{\beta}f \|_{L_{p}(\Gamma)}^{\alpha/\beta}
\|f\|_{L_{q}}^{1-\alpha/\beta},
\end{equation}
provided that
$$
\frac{1}{r}=\frac{\alpha}{\beta}\cdot\frac{1}{p}+
\Big(1-\frac{\alpha}{\beta}\Big)\frac{1}{q}.
$$
\end{lemma}

Proof. We basically mimic the proof
of Proposition 3.1 of \cite{Ad_75}.
Observe that \eqref{9.28.1} at the origin is easily obtained from summing up the inequalities in 
\eqref{1.17.2} and minimizing with respect 
to $\rho$. At any other point it is obtained
by changing the origin.
 Furthermore
clearly, we may assume that $f$ is bounded
with compact support.  Set $Q_{1}=\{(s,y):|y|\geq\sqrt s\}$,
$Q_{2}=\{(s,y):|y|\leq \sqrt s\}$.
Dealing with $P_{\alpha}(fI_{Q_{1}})$
we observe that $p_{\alpha}(s,r)\leq Nr^{-(d+2-\alpha)}$ if $r\geq\sqrt s$. Therefore,
$$
P_{\alpha}(fI_{Q_{1}\cap C^{c}_{\rho}})(0)
\leq N\int_{\rho}^{\infty}\frac{1}{r^{d+2-\alpha}}
\int_{0}^{r^{2}}\Big(\int_{|y|=r}f(s,y)\,d\sigma_{r}\Big)\,dsdr,
$$
where $d\sigma_{r}$ is the element of the surface area on $|y|=r$. By Lemma \ref{lemma 1.16.1}  ($\alpha<d+2$)
$$
P_{\alpha}(fI_{Q_{1}\cap C^{c}_{\rho}})(0)\leq N
\int_{\rho}^{\infty}\frac{1}{r^{d+3-\alpha}}
\int_{\rho}^{r}\Big(\int_{0}^{\rho^{2}}
\Big(\int_{|y|=\rho}f(s,y)\,d\sigma_{\rho}\Big)\,ds\Big)\,d\rho dr
$$
$$
\leq N
\int_{\rho}^{\infty}\frac{1}{r^{d+3-\alpha}}
\int_{0}^{r}\Big(\int_{0}^{r^{2}}
\Big(\int_{|y|=\rho}f(s,y)\,d\sigma_{\rho}\Big)\,ds\Big)\,d\rho dr
$$
$$
=N\int_{\rho}^{\infty}\frac{1}{r^{d+3-\alpha}}I(r)\,dr ,
$$
where  
$$
I(r)=\int_{C_{r}}f(s,y)\,dyds.
$$
We use that  $I(r)\leq Nr^{ d+2-\beta }M_{\beta}f(0)$ and that $\alpha<\beta$. Then we see that 
\begin{equation}
                          \label{1.17.3}
P_{\alpha}(fI_{Q_{1}\cap C^{c}_{\rho}})(0)\leq N\rho^{ \alpha-\beta }M_{\beta}f(0).
\end{equation}

Next, by using Lemma \ref{lemma 1.16.1}
we obtain that
$$
P_{\alpha}(fI_{Q_{2}\cap C^{c}_{\rho}})(0)
\leq
\int_{\rho^{2}}^{\infty}\frac{1}{s^{(d+2-\alpha)/2}}
\int_{|y|\leq\sqrt s}f(s,y)\,dyds
$$
$$
\leq N\int_{\rho^{2}}^{\infty}\frac{1}{s^{(d+4-\alpha)/2}}I(\sqrt s)\,ds=
N\int_{\rho}^{\infty}\frac{1}{r^{d+3-\alpha}}I(r)\,dr.
$$
This along with \eqref{1.17.3} prove  the second
inequality in \eqref{1.17.2}.

As long as the first inequality is concerned,
observe that  similarly to Lemma \ref{lemma 1.16.1} using that $\alpha>0$ we have
$$
P_{\alpha}(fI_{Q_{1}\cap C _{\rho}})(0)
\leq N\int_{0}^{\rho}\frac{1}{r^{d+2-\alpha}}
\int_{0}^{r^{2}}\Big(\int_{|y|=r}f(s,y)\,d\sigma_{r}\Big)\,dsdr
$$
$$
=N\int_{0}^{\rho}\frac{1}{r^{d+2-\alpha}}
\Big(\frac{\partial}{\partial r}\int_{0}^{r}
\Big(\int_{0}^{\tau^{2}}\int_{|y|=\tau}f(s,y)\,d\sigma_{\tau}\,ds\Big)\,d\tau\Big)dr
$$
$$
=J_{1}+ N
\int_{0}^{\rho}\frac{1}{r^{d+3-\alpha}}
\int_{0}^{r}\Big(\int_{0}^{\tau^{2}}
\Big(\int_{|y|=\tau}f(s,y)\,d\sigma_{\tau}\Big)\,ds\Big)\,d\tau dr
$$
$$
\leq J_{ 1}+N\int_{0}^{\rho}\frac{1}{r^{d+3-\alpha}}I(r)\,dr ,
$$
where
$$
J_{1}=N\frac{1}{\rho^{d+2-\alpha}}\int_{0}^{\rho}
\Big(\int_{0}^{\tau^{2}}\int_{|y|=\tau}f(s,y)\,d\sigma_{\tau}\,ds\Big)\,d\tau 
\leq N\frac{1}{\rho^{d+2-\alpha}}I(\rho)
$$
Here $I(r)
\leq Nr^{ d+2 }  Mf(0)$ and $\alpha>0$,
so that 
\begin{equation}
                          \label{1.17.4}
P_{\alpha}(fI_{Q_{1}\cap C _{\rho}})(0)
\leq N\rho^{\alpha }Mf(0).
\end{equation}
Furthermore,
$$
P_{\alpha}(fI_{Q_{2}\cap C _{\rho}})(0)
\leq N\int_{0}^{\rho^{2}}\frac{1}{s^{(d+2-\alpha)/2}}
\int_{|y|\leq\sqrt s}f(s,y)\,dyds
$$
$$
\leq J_{2}+ N\int_{0}^{\rho^{2}}\frac{1}{s^{(d+4-\alpha)/2}}I(\sqrt s)\,ds= J_{2}+
N\int_{0}^{\rho}\frac{1}{r^{d+3-\alpha}}I(r)\,dr,
$$
where
$$
J_{2}=N\frac{1}{\rho^{ d+2-\alpha }}\int_{0}^{\rho^{2}}\int_{|y|\leq \sqrt\tau}f(\tau,y)\,dyd\tau  \leq N\frac{1}{\rho^{d+2-\alpha}}I(\rho).
$$
This and \eqref{1.17.4} prove the first inequality in \eqref{1.17.2}.
 The lemma is proved.

\begin{remark}
                    \label{remark 1.18.1}
If $d=\alpha=1$ and $f$ is independent
of $t$, the inequalities \eqref{1.17.2}
and \eqref{9.28.1} are useless, because
the first one in \eqref{1.17.2} follows by
definition and the second one and 
\eqref{9.28.1} are trivial because $M_{\beta}f
=\infty$ ($\beta>\alpha=1$) unless $f=0$.

\end{remark}

If $\alpha$ is strictly less than the index of $L_{q}$, we have the following.
 \begin{corollary}
                    \label{corollary 9.29.1}
If $\alpha\in(0,(d+2)/q)$, $q\in(1,\infty)$,
  then
there exists a constant $N$
  such that for any
 $f\geq0$ we have
$$
\|P_{\alpha}f\|_{L_{r}}\leq N\|f\|_{L_{q}}
$$
as long as   
$$
\frac{d+2}{q}-\alpha=\frac{d+2}{r}.
$$
In particular, (a classical embedding result) if
  $1<q<d+2$ and $u\in C^{\infty}_{0}=
C^{\infty}_{0}(\bR^{d+1}) $, then
$$
\|Du\|_{L_{r}}\leq N\|\partial_{t}u+\Delta u\|_{L_{q}}\quad(\partial_{t}=\partial/\partial t)
$$
as long as
$$
\quad  \frac{d+2}{q}-1= \frac{d+2}{r}.
$$
\end{corollary}

Indeed, the first assertion follows from
H\"older's inequality and
\eqref{9.28.2} with $p=\infty$ and $\beta=(d+2)/q$ ($>\alpha$). The second assertion follows from the first one with $\alpha=1$ ($<\beta$) and the fact that
for $f=\partial_{t}u+\Delta u$ we have
$$
Du(t,x)=c\int_{\bR^{d +1 }_{+}}
\frac{y}{s^{(d+2)/2}}e^{-|y|^{2}/(4s)}
f(t+s,x+y)\,dyds,
$$
where $c$ is a constant  and
$(|y|/s^{1/2})e^{-|y|^{2}/(4s)}\leq
Ne^{-|y|^{2}/(8s)}$.

\begin{remark}
                          \label{remark 10.8.5}
After Corollary \ref{corollary 9.29.1} a natural
question arises as to what power of summability
$b=(b^{i})$ will be sufficient for
the term $b^{i}D_{i}u$ to be considered as a perturbation term in $\partial_{t}u+\Delta u
+b^{i}D_{i}u$ in the framework of the $L_{q}$-theory. Observe that, in the notation
of Corollary \ref{corollary 9.29.1} 
\begin{equation}
                            \label{1.29.3}
\|b^{i}D_{i}u\|_{L_{q}}\leq \|b\|_{L_{d+2}}
\|Du\|_{L_{r}}\leq N\|b\|_{L_{d+2}}\|\partial_{t}u+\Delta u\|_{L_{q}}.
\end{equation}
It follows that $b$ should be of class $L_{d+2}$ and $q<d+2$. Of course, if $b$
contains just bounded part, this part
in $b^{i}D_{i}u$ is taken care of by interpolation inequalities.

\end{remark}

In the next section we will also need
the following result.
\begin{corollary}
                     \label{corollary 1.17.1}
                           
For any $\alpha\in(0,\beta ),\beta\in(0,d+2]$ there exists a constant $N$
  such that for any $g\geq0$, $\rho\in(0,\infty)$, and  $(t,x)\in C_{\rho}$ we have
$$
P_{\alpha}(I_{C^{c}_{2\rho}}g)(t,x)
\leq N\rho^{\alpha-\beta}M_{\beta}
g(t,x).
$$

\end{corollary}

Indeed, since 
$$
\{t+s\geq 4\rho^{2}\quad\text{or}
\quad |x+y|\geq 2\rho\}\subset
\{ s\geq  \rho^{2}\quad\text{or}
\quad | y|\geq  \rho\}
$$
for $f=g(t+\cdot,x+\cdot)$ we have 
$$
P_{\alpha}(I_{C^{c}_{2\rho}} g)(t,x)
\leq\int_{\bR^{d+1} }I_{C^{c}_{\rho}}
(s,y)p_{\alpha}(s,y)g(t+s,x+y)\,dyds=
P_{\alpha}(I_{C^{c}_{\rho}}f)(0) 
$$
$$
\leq N\rho^{\alpha-\beta}M_{\beta}
f(0)=N\rho^{\alpha-\beta}M_{\beta}
g(t,x).
$$

\mysection{A parabolic analog of the Adams
Theorem 3.1 of \protect\cite{Ad_75}}

                       \label{section 10.8.2}

\begin{theorem}
                       \label{theorem 9.29.1}
For any $\alpha\in(0,\beta ),\beta\in(0,(d+2)/q]$, 
$q\in(1,\infty)$, and $r$ such that
$$
r(\beta-\alpha)=q\beta ,
$$
there is a constant $N$ such that for any
   $f\geq0$ we have
\begin{equation}
                          \label{9.29.3}
 \|P_{\alpha}f\|_{E_{r,\beta-\alpha}}
\leq N \|f\|_{E_{q,\beta}}.
\end{equation}
\end{theorem}

Proof. It suffices to prove that for any $\rho>0$
$$
\rho^{\beta-\alpha}\Big(\dashint_{C_{\rho}}|P_{\alpha}f|^{r}\,dz\Big)^{1/r}
\leq N\|f\|_{E_{q,\beta}},
$$
that is
\begin{equation}
                          \label{9.29.4}
\rho^{\beta-\alpha-(d+2)/r}\Big(\int_{C_{\rho}}|P_{\alpha}f|^{r}\,dz\Big)^{1/r}
\leq N\|f\|_{E_{q,\beta}},
\end{equation}

Observe that by H\"older's inequality
$M_{\beta}f\leq N\|f\|_{E_{q,\beta}}$ and by definition
$$
\Big(\int_{\bR^{d+1}}I_{C_{2\rho}}f^{q}\,dz
\Big)^{1/q}
\leq N\rho^{(d+2)/q-\beta  }\|f\|_{E_{q,\beta}}.
$$
It follows from  Lemma \ref{lemma 9.28.1} 
with $p=\infty$ that
$$
\Big(\int_{C_{\rho}}|P_{\alpha}(I_{C_{2\rho}}f)|^{r}
dz\Big)^{1/r}\leq N\rho^{((d+2)/q-\beta)(1-\alpha/\beta) }\|f\|_{E_{q,\beta}}
$$
$$
=N\rho^{(d+2)/r-\beta+\alpha}\|f\|_{E_{q,\beta}}.
$$
Furthermore, by Corollary \ref{corollary 1.17.1}
$$
\Big(\int_{C_{\rho}}|P_{\alpha}(I_{C^{c}_{2\rho}}f)|^{r}
dz\Big)^{1/r}\leq N\rho^{(d+2)/r} \sup_{C_{\rho}}P_{\alpha}(I_{C^{c}_{2\rho}}f)
$$
$$
\leq N\rho^{(d+2)/r+\alpha-\beta}E_{q,\beta}f.
$$
By combining these estimates we come to 
\eqref{9.29.4} and the theorem is proved.

\begin{remark}
                    \label{remark 9.29.2}
We did not explicitly used that $\beta
\leq (d+2)/q$ and formally the proof is valid
for any $\beta\in(0,\infty)$
  if in Definition \ref{8.11.2}
we allow any $\beta>0$.  However, if $\beta>(d+2)/q$
and $f\ne 0$, the right-hand side of \eqref{9.29.3} is infinite. Therefore,
to make Theorem \ref{theorem 9.29.1}
nontrivial one requires $\beta\leq (d+2)/q$.
\end{remark}

\begin{remark}
There is a simple relation of $P_{\alpha_{1}}
P_{\alpha_{2}}$ to $P_{\alpha_{1}+\alpha_{2}}$,
which, in light of Theorem \ref{theorem 9.29.1}, implies that, if $\beta>\alpha_{2}\geq\alpha_{1}>0$,
$q_{1},q_{2}\in(1,\infty)$, $q_{1}(\beta-\alpha_{1})=q_{2}(\beta-\alpha_{2})
\leq d+2$, then
$\|P_{\alpha_{2}}f\|_{E_{q_{2},\beta-\alpha_{2}}}
\leq N\|P_{\alpha_{1}}f\|_{E_{q_{1},\beta-\alpha_{1}}}$.
We leave  details of the proof to the reader and
we do not use this fact in what follows.
\end{remark}
 
The following, obtained similarly to Corollary \ref{corollary 9.29.1}, was communicated
to the author by Hongjie Dong.
\begin{corollary}
                 \label{corollary 9.30.1}
If   $1<q<d+2$, $\beta\in(1,(d+2)/q]$, and $u\in C^{\infty}_{0} $, then
$$
\|Du\|_{E_{r,\beta-1}}\leq N\|\partial_{t}u+\Delta u\|_{E_{q,\beta}}
$$
as long as  
\begin{equation}
                                  \label{10.2.1}
r(\beta-1)=q\beta,\quad\text{\rm that is}\quad
\frac{1}{r}=\frac{1}{q}-\frac{1}{\beta q}.
\end{equation}
\end{corollary}

\begin{remark}
                      \label{remark 10.9.1}
For $\beta=(d+2)/q$ Corollary \ref{corollary 9.30.1} yields 
the second part of Corollary \ref{corollary 9.29.1} once more. This is because
$ E_{q,(d+2)/q} = L_{q} $.
\end{remark}

\begin{remark}
                      \label{remark 10.8.3}
In the framework of the Morrey spaces
Corollary \ref{corollary 9.30.1} opens up the
possibility to treat the terms like
$b^{i}D_{i}u$ as perturbation terms
in operators like $\partial_{t}u+\Delta u
+b^{i}D_{i}u$ even with rather low
summability properties of $b=(b^{i})$.
To show this, observe that for $f,g\geq0$
in the notation of Corollary \ref{corollary 9.30.1}
$$
\rho^{\beta}\dashnorm fgI_{C_{\rho}}\|_{L_{q}}\leq \rho
\dashnorm f I_{C_{\rho}}\|_{L_{\beta q}}\cdot \rho^{\beta-1}\dashnorm gI_{C_{\rho}}\|_{L_{r}}.
$$
It follows that
\begin{equation}
                                \label{1.29.2}
\|b^{i}D_{i}u\|_{E_{q,\beta}}
\leq \|b \|_{E_{\beta q,1}}\| D u\|_{E_{r,\beta-1}}
\leq N\|b \|_{E_{\beta q,1}}\|\partial_{t}u+\Delta u\|_{E_{q,\beta}}.
\end{equation}
For $\beta=(d+2)/q$ estimate \eqref{1.29.2}
coincide with \eqref{1.29.3}, but for $\beta
<(d+2)/q$ in the framework of Morrey spaces we allow
  $b$ to be summable to the power
$\beta q<d+2$ in contrast with Remark 
\ref{remark 10.8.5}.
  However, we need $\|b \|_{E_{\beta q,1}}<\infty$ and, if we ask ourselves what $r$
should be in order for $b\in L_{r}$ to
have $\|b \|_{E_{\beta  q,1}}<\infty$, the answer is
$r=d+2$ at least. Still we gain the possibility to have higher singularities of $b$ than functions from $L_{d+2}$. Elliptic versions of
\eqref{1.29.2}  for usual or generalized Morrey spaces are found in many papers,
see, for instance, \cite{FHS_17} and the references therein.

\end{remark}

Next we move to deriving ``local'' versions of the above results.
A statement somewhat weaker than Corollary \ref{corollary 9.30.1} can be obtained from the following   general result
by taking $(S,T)$ to be large enough
and then sending $S\to-\infty,T\to\infty$.

\begin{theorem}
                        \label{theorem 10.2.1}
Let   
$1<q<d+2$, $\beta\in(1,(d+2)/q]$
and let \eqref{10.2.1} hold.   Then
there is a constant $N$ such that
for any $u\in C^{\infty}_{0}$, $-\infty<S<T<\infty$, and 
$Q_{S,T}=(S,T)\times\bR^{d}$
\begin{equation}
                                  \label{10.2.2}
 \|Du\|_{E_{r,\beta-1}(Q_{S,T})}\leq N\||\partial_{t}u|+|\Delta u|\,\|_{E_{q,\beta}(Q_{S,T})}
+N(T-S)^{-1}\|u\|_{E_{q,\beta}(Q_{S,T})}.
\end{equation}

\end{theorem}

Proof. Shifting and changing the scales in $\bR^{d+1}$ allow  us to assume that $S=-1=-T$.
 In that case consider the mapping $\Phi:[-3/2, 3/2]\to[-1,1]$, $\Phi(t)=t\big(2/(|t|\vee 1)-1\big)$ that   preserves $[-1,1]$,  
is Lipschitz continuous
and has Lipschitz continuous inverse if restricted to $[-3/2, 3/2]\setminus (-1,1)$.
  Then, obviously,  for $w(t,x)=v(\Phi (t),x)$
we have
\begin{equation}
                                \label{5.22.4}
\|w I_{Q_{-3/2,3/2}}
\|_{E_{q,\beta} }
\leq N \|v 
\|_{E _{q,\beta}(Q_{-1,1})},
\end{equation}
where $N=N( q)$.
 
Now take $(t,x)\in Q_{-1,1}$, $\rho\in(0,\infty)$,
and take $\zeta\in C^{\infty}_{0}(\bR )$
such that $\zeta=1$ on $(-1,1)$, $\zeta=0$
outside $(-3/2,3/2)$, and $|\zeta|+| \zeta'| 
\leq 4 $.

Although the function $\zeta w$,
where $w(s,y)=u(\Phi(s),y)$, is not as smooth as
  required in Corollary \ref{corollary 9.30.1}
the argument leading to it applies
to $\zeta(s)w(s,y)$ (we have a general Remark
\ref{remark 10.10.4} to that effect) and since $r(\beta-1)=q\beta$ we have
$$
\rho^{\beta-1}\dashnorm DuI_{Q_{-1,1}}\|_{L_{r}( C_{  \rho}(t,x))}
\leq N\rho^{\beta-1}\dashnorm D(\zeta w)\|_{L_{r}( C_{ \rho}(t,x))}
$$
$$
\leq N\|I_{Q_{-3/2,3/2}}(|\partial_{t} (\zeta w)|
+|\zeta\Delta w|)\|_{E_{q,\beta}}.
$$
It only remains to note that the last expression
is less than the right-hand side of 
\eqref{10.2.2}
in light of \eqref{5.22.4}. The theorem is proved.

To prove an interpolation theorem in $C_{R}$
we need two lemmas.
\begin{lemma}
                    \label{lemma 10.4.1}
Let $0< R_{1}<1<R_{2}< \infty $, 
$1\leq q <\infty$, $\beta\in(0,(d+2)/q]$. Define
$\Gamma_{1}=\bar B_{1}\setminus B_{R_{1}}$, 
$\Gamma_{2}=\bar B_{R_{2}}\setminus B_{1}$ and let  
  $\Phi:\Gamma_{2}
\to \Gamma_{1}$ be a smooth
one-to-one mapping with $|D\Phi|,|D\Phi^{-1}|\leq K$, where $K$ is a constant. Let $v(t,x)\geq0$
be zero outside $G_{2}:=(0,1)\times\Gamma_{2}$
and set $u(t,x)=v(t,\Phi^{-1}(x))I_{\Gamma_{1}}(x)$. Then
\begin{equation}
                                \label{10.4.1}
\|v\|_{E_{q,\beta}((0,1)\times B_{R_{2}})}
\leq N(d,q,\beta,K)\|u\|_{E_{q,\beta}(C_{1})}.
\end{equation}
\end{lemma}

Proof. Take $(t,x)\in (0,1)\times B_{R_{2}}$
and $\rho>0$. Then  
$$
\rho^{\beta}\Big(\frac{1}{\rho^{d+2}}\int_{C_{\rho}(t,x)}
I_{(0,1)\times B_{R_{2}}}v^{q}\,dyds\Big)^{1/q}
$$
$$
\leq N\rho^{\beta}\Big(\frac{1}{\rho^{d+2}}
\int_{\Psi(C_{\rho}(t,x)\cap G_{2})}
I_{C_{1}}u^{q}\,dyds\Big)^{1/q}=:I,
$$
where $\Psi(s,y)=(s,\Phi(y))$.
Observe that, if $C_{\rho}(t,x)\cap G_{2}
\ne\emptyset$, then $|y_{1}-y_{2}|\leq 2\rho$
for any $y_{1},y_{2}\in
C_{\rho}(t,x)\cap G_{2}$. It follows that
$\Phi(C_{\rho}(t,x)\cap G_{2})\subset B$,
where $B$ is a ball
of radius $2K\rho$ with center in $B_{1}$,
 and
$$
I\leq N(2K\rho)^{\beta}\Big(\frac{1}{(2K\rho)^{d+2}}
\int_{(t,t+(2K\rho)^{2})
\times B}I_{C_{1}}u^{q}\,dyds\Big)^{1/q}
\leq N\|u\|_{E_{q,\beta}(C_{1})}.
$$
This proves the lemma.

The following lemma about the interpolation
inequality \eqref{8.4.3} is quite natural
and obviously useful,
but its elliptic counterpart was proved
only rather late in \cite{Kr_303}.
One of its goals is to be able to treat
$b^{i}D_{i}u$, when $b$ is bounded,
as a perturbation term.

\begin{lemma}
                      \label{lemma 7.30.1}
Let   $p\in(1,\infty)$, $0<\beta\leq (d+2)/p$. Then there is a constant
$N $ such that, for any 
$R\in(0,\infty)$, $\rho\leq 2R  $,
$C\in \cC_{\rho}$ with its geometric center
in $C_{R}$,
    $\varepsilon
\in(0,1] $, 
and $u\in C^{\infty}_{0}$, we have
$$
 \rho^{\beta}
\dashnorm I_{C_{R}}Du\|_{L_{p}(C )}
\leq 
 N\varepsilon R \sup_{\rho\leq s\leq 2R}s^{\beta}
\dashnorm I_{C_{ R}}(|\partial_{t}u|+|D^{2}u|)\|_{L_{p}(C(s))}
$$
\begin{equation}
                              \label{7.30.1}
+N\varepsilon^{-1}R^{-1}
\sup_{\rho\leq s\leq 2R}s^{\beta}\dashnorm I_{C_{R}}( u-c)\|_{L_{p}(C(s))},
\end{equation}
where $c$ is any constant and
$C(s)\in \cC_{s}$ with the geometric center
the same as $C$. In particular,
\begin{equation}
                              \label{8.4.3}
\|Du\|_{E_{p,\beta}(C_{R})}\leq
N\varepsilon R \||\partial_{t}u|+|D^{2}u|\,\|_{E_{p,\beta}(C_{R})}
+N\varepsilon^{-1}R^{-1}
\| u\|_{E_{p,\beta}(C_{R})}.
\end{equation}
\end{lemma}

Proof. Changing scales shows that we may assume that $R=1$. Obviously we may also assume that $c=0$. Then denote $v= Du$,
$w= |\partial_{t}u|+|D^{2}u| $, $G_{s}=C(s) \cap C_{1}$,
$$
U=\sup_{\rho\leq s\leq 2 }s^{\beta}
\dashnorm  u\|_{L_{p}(G_{s})},\quad
W=\sup_{\rho\leq s\leq 2 }s^{\beta}
\dashnorm  (|\partial_{t}u|+|D^{2}u|)\|_{L_{p}(G_{s})},
$$
 
By Poincar\'e's inequality
(see, for instance, Lemma \ref{lemma 10.10.1}), for $ \rho\leq s\leq  2$,
$$
\dashnorm v-v_{G_{s}}\|_{L_{p}(G_{s})}\leq
N(d,p)s\dashnorm w\|_{L_{p}(G_{s})}
\leq Ns^{1-\beta}W.
$$
Also by interpolation inequalities,
   there
exists a constant $N=N( d,p)$
such that, for  
 $\varepsilon \in(0,1]$ and $\varepsilon\leq s\leq  2$ ,
$$
\dashnorm v-v_{G_{s}}\|_{L_{p}(G_{s})}
\leq 2\dashnorm v\|_{L_{p}(G_{s})}
\leq N \dashnorm w\|^{1/2}_{L_{p}(G_{s})}
\dashnorm  u\|^{1/2}_{L_{p}(G_{s})}
$$
\begin{equation}
                                 \label{7.31.1}
+Ns^{-1}\dashnorm  u\| _{L_{p}(G_{s})}
\leq N \dashnorm w\|^{1/2}_{L_{p}(G_{s})}
\dashnorm  u\|^{1/2}_{L_{p}(G_{s})}
+N\varepsilon^{-1}\dashnorm  u\| _{L_{p}(G_{s})},
\end{equation}
which for $2\geq s\geq \varepsilon\vee \rho$    yields
$$
s^{\beta}\dashnorm v-v_{G_{s}}\|_{L_{p}(G_{s})}
\leq N W^{1/2}U^{1/2}   
+N \varepsilon^{-1} U.
$$
  Hence, for any $\varepsilon\in(0,1]$ 
and $\rho\leq s\leq 2$
$$
s^{\beta}\dashnorm v-v_{G_{s}}\|_{L_{p}(G_{s})}
\leq N_{1}\varepsilon W 
+N_{2}\varepsilon^{-1}U,
$$
where $N_{1}=N_{1}(d,p)$, 
$N_{2}=N_{2}( d,p)$.

Following Campanato,
one can transform this result
to estimate $v_{G_{s}}$ going along $\rho$,
$2\rho$,... and,
since $\beta\in(0,(d+2)/p]$, by Campanato's
results (cf.~for instance, Proposition 5.4 in \cite{MM_12}) one gets  that  
$$
\rho^{\beta}\dashnorm v \|_{L_{p}(G_{\rho})}
\leq N_{3}( N_{1}\varepsilon W 
+N_{2}\varepsilon^{-1}U)+N_{3}\dashnorm v \|_{L_{p}(G_{2})},
$$
where $N_{3}=N_{3}(d,p,\beta)$. 
We estimate the last term as in \eqref{7.31.1}
and come to what implies
\eqref{7.30.1}.
The lemma is proved.

The following is a local version
of Corollary \ref{corollary 9.30.1}.
It allows us to draw the same conclusions as in Remark \ref{remark 10.8.3}
in bounded domains.

\begin{theorem}
                        \label{theorem 10.2.2}
Let  
$1<q<d+2$, $\beta\in(1,(d+2)/q]$
and let $r(\beta-1)=q\beta$.  
 Then  there is a constant $N$ such that
for any $R\in(0,\infty]$, $u\in C^{\infty}_{0}$,
\begin{equation}
                                  \label{10.2.20}
 \|Du\|_{E_{r,\beta-1}(C_{R})}\leq N\||\partial_{t}u|+|D^{2} u|\,\|_{E_{q,\beta}(C_{R})}
+NR^{-2}\|u\|_{E_{q,\beta}(C_{R})}.
\end{equation}

\end{theorem}

Proof. The case of $R=\infty$ is obtained by passing to the limit. In case $R<\infty$,
as usual, we may assume that $R=1$.
In that case, mimicking the Hestenes formula,
for $1\leq |x|\leq 6/5$ define
$$
v(t,x)=6u(t,x(2/|x|-1))
-8u(t,x(3/|x|-2))+3u(t,x(4/|x|-3))
$$
$$
=:6v_{1}-8v_{2}
+3v_{3}
$$
and for $|x|\leq 1$ set $v(t,x)=u(t,x)$.
One can easily check that $v\in C^{1,2}([0,1]
\times B_{6/5})$. In light of Lemmas
\ref{lemma 10.4.1} and \ref{lemma 7.30.1}, for instance,
$$
\|D^{2} v\|_{E_{q,\beta}((0,1)
\times B_{6/5})}\leq \|D^{2} u\|_{E_{q,\beta}(C_{1})}+N\|I_{B_{6/5}\setminus B_{1}}D^{2} v_{1}\|_{E_{q,\beta}((0,1)
\times B_{6/5})}+...
$$
$$
+N\|I_{B_{6/5}\setminus B_{1}}D^{2} v_{3}\|_{E_{q,\beta}((0,1)
\times B_{6/5})}\leq
N\|D^{2} u\|_{E_{q,\beta}(C_{1})}
+N\|Du\|_{E_{q,\beta}(C_{1})} 
$$
\begin{equation}
                                \label{10.4.3}
\leq
N\|D^{2} u\|_{E_{q,\beta}(C_{1})}
+N\| u\|_{E_{q,\beta}(C_{1})}.
\end{equation}
  
Now take $(t,x) \in C_{1}$, $\rho\in(0,\infty)$,
and take $\zeta\in C^{\infty}_{0}(\bR^{d})$
such that $\zeta=1$ on $B_{1}$, $\zeta=0$
outside $B_{6/5}$, and $|\zeta|+|D\zeta| 
+|D^{2}\zeta|\leq N=N(d)$.

By using  Theorem \ref{theorem 10.2.1} we get
$$
\rho^{\beta-1}\dashnorm DuI_{C_{1}}\|_{L_{r}( C_{  \rho}(t,x ) )}
\leq N\rho^{\beta-1}\dashnorm I_{Q_{0,1}}D(\zeta v) \|_{L_{r}( C_{ \rho}(t,x))}
$$
$$
\leq N\|D (\zeta v)\|_{E_{r,\beta-1}(Q_{0,1})}
\leq N\|\,|\partial_{t}(\zeta v)|+|\Delta(\zeta v)|\,\|_{E _{q,\beta}( Q_{0,1} )}
+N\|\zeta v\|_{E _{q,\beta}( Q_{0,1} )}
$$
$$
\leq  N\|\,|\partial_{t}(\zeta v)|+|\Delta(\zeta v)|\,\|_{E_{q,\beta}((0,1)
\times B_{6/5})}
+N\|v\|_{E_{q,\beta}((0,1)
\times B_{6/5})}
$$

It only remains to note that the last expression
is less than the right-hand side of 
\eqref{10.2.20} as is well seen from
  \eqref{10.4.3}. The theorem is proved.

\begin{remark}
                       \label{remark 2.16.1}
By considering functions depending only on
$x$ we naturally obtain ``elliptic'' analogs
of our results. For instance, for $G\subset
\bR^{d}$
by defining
$$
\|g\|_{E_{p,\beta}(G)}=\sup_{\rho<\infty,
x\in G}\rho^{\beta}\dashnorm gI_{G}\|_{L_{p}(B_{\rho}(x))},
$$
we get from \eqref{10.2.20} for $u\in C^{\infty}_{0}(\bR^{d})$ that
\begin{equation} 
                       \label{2.16.3}
 \|Du\|_{E_{r,\beta-1}(B_{R})}\leq N\||  D^{2} u\,\|_{E_{q,\beta}(B_{R})}
+NR^{-2}\|u\|_{E_{q,\beta}(B_{R})},
\end{equation}
whenever $1<q<d,\beta\in(1,d/q]$ and $r(\beta-1)=q\beta$. Actually, formally,
one gets \eqref{2.16.3} even for $\beta\leq
(d+2)/q$, but for $\beta >d/q$, both sides
of \eqref{2.16.3} are infinite unless $u=0$.

After that arguing as in \eqref{1.29.2} we see
that for $1<q<d,\beta\in(1,d/q]$
\begin{equation}
                                \label{2.16.4}
\|b^{i}D_{i}u\|_{E_{q,\beta}(B_{1})}
\leq N\|b \|_{E_{\beta q,1}(B_{1})}\| \Delta u\|_{E_{q,\beta}(B_{1})}+N
\|u\|_{E_{q,\beta}(B_{1})}.
\end{equation}

From the point of view of the theory of elliptic equations
the most desirable version of \eqref{2.16.4}
would be 
\begin{equation}
                            \label{2.15.1}
\|b^{i}D_{i}u\|_{E_{q,\beta}(B_{1})}
\leq \varepsilon\| \Delta u\|_{E_{q,\beta}(B_{1})}+N(\varepsilon)
\|u\|_{E_{q,\beta}(B_{1})}
\end{equation}
for any $\varepsilon>0$ with   $N(\varepsilon)$ independent
of $u$. This fact is, actually, claimed
in Theorem 5.4 of \cite{FHS_17}.
We will show that \eqref{2.15.1}
cannot hold if $\varepsilon$ is small enough. 

Let $h(t)$ be a smooth nondecreasing function on $\bR$
such that $h(t)=0$ for $t\leq 0$, $h(t)=t$
for $t\geq 1$ and for $\delta>0$ set
$u_{\delta}(x)=h(\ln(\delta/|x|))$. Let
$1<q<d/2$,   $\beta=2$,   $b(x)=1/|x|$.

Then
$$
\|u\|_{E_{q,\beta}(B_{1})}
\leq N(d)\|u_{\delta}\|_{L_{d}(B_{1})}\to0
$$
as $\delta\downarrow0$. At the same time
$$
D_{i}u_{\delta}=-\frac{x_{i}}{|x|^{2}}h',\quad
D_{ij}u_{\delta}=\frac{1}{|x|^{2}}
\Big(2\frac{x_{i}x_{j}}{|x|^{2}}-\delta_{ij} 
\Big)h'+\frac{1}{|x|^{2}}\frac{x_{i}x_{j}}{|x|^{2}}h''.
$$
It is seen that $|D^{2}u_{\delta}|\leq N(d)/|x|^{2}$
and, since $q<d/2$, the $E_{q,\beta}(B_{1})$-norm of $D^{2}u_{\delta}$ is bounded as $\delta
\downarrow 0$. Also, for $|x|\leq\delta/e$,
we have $b|Du_{\delta}|=1/|x|^{2}$, so that
for $r\leq\delta/e$
$$
\Big( \dashint_{|x|\leq r}b^{q}|Du_{\delta}|^{q}\,dx\Big)^{1/q}=N(d,p)r^{-2}.
$$
It follows that the $E_{q,\beta}(B_{1})$-norm of $b|Du_{\delta}|$ is bounded away from zero as $\delta
\downarrow 0$ and this shows that \eqref{2.15.1} cannot hold for all $\delta>0$
if $\varepsilon $ is small enough.
\end{remark}

\mysection{A parabolic version
of Chiarenza--Frasca result \cite{CF_90}}

                       \label{section 1.25.1}

In Remark \ref{remark 10.8.3} we have shown
how to estimate a Morrey norm of $|b|\,|Du|$
in terms of a Morrey norm of $b$.
Here, following \cite{CF_90}, we show
how to estimate an $L_{p}$-norm of the same quantity through the $L_{p}$-norms of $\partial_{t}u$ and $D^{2} u$.  
\begin{theorem}
                    \label{theorem 1.18.1}
Let $d+2\geq q> p>1$,   $b\in E_{q,1}$.  Then for any $f\geq0$
we have
\begin{equation}
                         \label{5.30.1}
I:=\int_{\bR^{d+1}}|b|^{p}(P_{1}f)^{p}\,dz\leq
N\|b\|_{E_{q,1}}^{p}\|f\|_{L_{p}},
\end{equation}
where $N$ depends only on $d,p,q$.
In particular  (see the proof of
Corollary \ref{corollary 9.29.1}), 
 for any $u\in C^{\infty}_{0}$
\begin{equation}
                       \label{1.18.3}
 \int_{\bR^{d+1} }|b |^{p}|Du |^{p}\,dz
\leq N\|b\|_{E_{q,1}}^{p}K,
\end{equation}
where
$
K=\|D^{2}u,\partial_{t}u\|^{p}_{L_{p}}
$
and $N$ depends only on $d,p,q$.

\end{theorem}

Observe that we already know this result if
$q=d+2$ from Remarks \ref{remark 10.8.5}
or \ref{remark 10.8.3}. 

In the proof we are going to use ``parabolic" versions of some results from Real Analysis
associated with balls  and cubes. These versions are obtained
by easy adaptation  of the corresponding arguments by replacing balls with
parabolic  cylinders and cubes with parabolic boxes. To make the adaptation more natural
we introduce the ``symmetric'' maximal
parabolic function
operator by
$$
\hat Mf(t,x)=\sup_{\substack{C\in\cC,\\   C\ni (t,x)}}\dashint_{C}|f|\,dz,
$$
where (recall that) $\cC$ is the set of $C_{r}(z)$, $r>0$,
$z\in \bR^{d+1}$.
To prove the theorem we need the following.

\begin{lemma}
                        \label{lemma 1.22.1}
a) For  $r\in(0,\infty)$ define 
$D _{r}
=\{|t|\leq r^{2},|x|\leq r \}$. Then
\begin{equation}
                             \label{1.22.2}
\hat MI_{D_{r}}(t,x)\leq I_{D_{2r}}+N
I_{D^{c}_{2r}} \frac{r^{d+2}}{  |t| ^{(d+2)/2}
\vee |x| ^{d+2}}\leq N^{2}\hat MI_{D_{r}}(t,x),
\end{equation}
where $N=N(d)$.

b) For any nonnegative $g(t,x)$, $q\in[1,\infty)$, $\beta\in(0,d+2]$,
$\alpha> 0$, $\alpha>1-q\beta/(d+2)$, and $r\in(0,\infty)$
\begin{equation}
                             \label{1.22.3}
\int_{\bR^{d}} g ^{q}\big(\hat MI_{D_{r}}
\big)^{\alpha}\,dz
\leq N(d,q,\alpha,\beta)r^{d+2-q\beta}\|g\|^{q}_{E_{q,\beta}}.
\end{equation}
\end{lemma}

Proof. Assertion a) is proved by elementary
means. To prove b), we use a) and split $D^{c}_{2r} $
into two parts $D^{c}_{2r}\cap\{|x|^{2}\geq |t|\} $ and $D^{c}_{2r}\cap\{|x|^{2}< |t|\} $
and, taking into account obvious symmetries, we see that it suffices to show that
$$
I_{1}:=\int_{4r^{2}}^{\infty}\int_{B_{\sqrt{ t}}}\frac{g^{q}(t,x)}{ t ^{\alpha(d+2)/2}}\,dxdt
\leq Nr^{(d+2)(1-\alpha)-q\beta}\|g\|^{q}_{E_{q,\beta}},
$$
$$
I_{2}:=\int_{|x|\geq 2r }\int_{ 0}^{|x|^{2}}\frac{g^{q}(t,x)}{|x|^{\alpha(d+2)}}\,dtdx
\leq Nr^{(d+2)(1-\alpha)-q\beta}\|g\|^{q}_{E_{q,\beta}}.
$$
 
By observing that
$$
\frac{1}{t^{\alpha(d+2)/2}}\int_{4r^{2}}^{t}\Big(\int_{B_{\sqrt s}}  g^{q}(s,x)\,dx\Big)ds\leq \frac{t^{(d+2)/2-q\beta/2}}{t^{\alpha(d+2)/2}}\|g\|^{q}_{E_{q,\beta}}\to0
$$
as $t\to\infty$, we have
$$
I_{1}=\int_{4r^{2}}^{\infty}\frac{1}{t^{\alpha(d+2)/2}}
\frac{d}{dt}\Big(\int_{4r^{2}}^{t}\Big(\int_{B_{\sqrt s}}  g^{q}(s,x)\,dx\Big)ds
\Big)\,dt
$$
$$
= N\int_{4r^{2}}^{\infty}\frac{1}{t^{\alpha(d+2)/2+1}}\Big(\int_{4r^{2}}^{t}\int_{B_{\sqrt s}}  g^{q}(s,x)\,dxds
\Big)\,dt
$$
$$
\leq N\|g\|^{q}_{E_{q,\beta}}\int_{4r^{2}}^{\infty}
\frac{t^{(d+2)/2-q\beta/2}}{t^{\alpha(d+2)/2+1}}\,dt
=Nr^{(d+2)(1-\alpha)-q\beta}\|g\|^{q}_{E_{q,\beta}}.
$$

Also as is easy to see  
$$
I_{2}=N\int_{2r }^{\infty}\frac{1}{\rho^{\alpha(d+2)}}
\int_{0}^{ \rho^{2}}\Big(\int_{|x|=\rho}
  g^{p}(t,x)\,d\sigma_{\rho}\Big)\,dtd\rho
$$
$$
\leq N\int_{2r }^{\infty}\frac{1}{\rho^{\alpha(d+2)}}\frac{\partial}{\partial\rho}\Big(
\int_{4r^{2} }^{ \rho^{2}}\Big(\int_{|x|\leq\rho}
  g^{p}(t,x)\,dx\Big)\,dt\Big)d\rho
$$
$$
\leq N\|g\|^{q}_{E_{q,\beta}}\int_{2r }^{\infty}\frac{\rho^{d+2-q\beta}}{\rho^{\alpha(d+2)+1}}\,d\rho
=Nr^{(d+2)(1-\alpha)-q\beta}\|g\|^{q}_{E_{q,\beta}}.
$$
This proves the lemma.

{\bf Proof of Theorem \ref{theorem 1.18.1}}. We follow  some  arguments in \cite{CF_90} and may assume that $b\geq0$.
First  set $r_{0}=(p+q)/2$ and  assume that there is a constant $N_{0}$ such that $\hat M(|b|^{r_{0}})\leq N_{0}|b|^{r_{0}}$, that is,   $|b|^{r_{0}}$ is in the class $A_{1}$
of Muckenhoupt.
Observe that by H\"older's inequality $\|b\|_{E_{r_{0},1}}\leq \|b\|_{E_{q,1}}$. 
It is convenient to prove the following version 
of \eqref{5.30.1} (notice $r_{0}$ in place of $q$)

\begin{equation}
                         \label{5.30.6}
I \leq
N\|b\|_{E_{r_{0},1}}^{p}\|f\|_{L_{p}},
\end{equation}

Then assume that $b\geq 0$, set
$u=P_{1}f$, and write   
\begin{equation}
                           \label{5.30.2}
I=\int_{\bR^{d+2}}\big(b^{p}u^{p-1}\big) P_{1}f
\,dz=\int_{\bR^{d+2}}P_{1}^{*}\big(b^{p}u^{p-1}\big)  f
\,dz\leq \|f\|_{L_{p}}\big\|
P_{1}^{*}\big(b^{p}u^{p-1}\big)\big\|_{L_{p'}},
\end{equation}
where $p'=p/(p-1)$ and $P^{*}_{1}$ is the conjugate operator for $P_{1}$, namely, for any $g\geq0$,   
\begin{equation}
                           \label{5.30.1}
(P^{*}_{1}g)
(s,x)=\big(P_{1}(g(-\cdot,-\cdot)\big)(-s,-x).
\end{equation}

Next, take $\gamma>0$, such that $(1+\gamma)p
\leq r_{0}$, $1+\gamma p'\leq r_{0}$,
 and $p\geq 1+\gamma$. Note that
$$
P_{1}^{*}\big(b^{p}u^{p-1}\big)=
P_{1}^{*}\big(b^{1+\gamma}\big(b^{p-1-\gamma}
u^{p-1}\big)\big)
$$
$$
\leq\Big(P_{1}^{*}\big(b^{(1+\gamma)p})\Big)^{1/p}
\Big( P_{1}^{*}\big(b^{p-\gamma p'}u^{p}\Big)^{(p-1)/p}.
$$
It follows that 
$$
\big\|
P_{1}^{*}\big(b^{p}u^{p-1}\big)\big\|_{L_{p'}}
\leq\Big(\int_{\bR^{d}}b^{p-\gamma p'}u^{p}
P_{1}\Big[\Big(P_{1}^{*}\big(b^{(1+\gamma)p})\Big)^{1/(p-1)}\Big]\,dz\Big)^{(p-1)/p}.
$$
Now in light of \eqref{5.30.2} we see that, to prove \eqref{5.30.6} in our particular case, it only remains to show that
\begin{equation}
                             \label{5.30.4}
P_{1}\Big[\Big(P_{1}^{*}\big(b^{(1+\gamma)p})\Big)^{1/(p-1)}\Big]\leq Nb^{\gamma p'}\|b\|^{p'}_{E_{r_{0},1}}.
\end{equation}
For $\alpha=1$ and $\beta=(1+\gamma)p$ ($>\alpha$) it follows from \eqref{9.28.1}
and \eqref{5.30.1} that 
$$
P_{1}^{*}\big(b^{(1+\gamma)p}\big)\leq N\|b\|_{E_{\beta,1}}\big(\hat M\big(b^{(1+\gamma)p}\big)\big)^{1-1/(p+\gamma p)}.
$$
where the last factor by assumption 
(and H\"older's inequality) is dominated by $Nb^{(1+\gamma)p-1}$
and $\|b\|_{E_{\beta,1}}\leq \|b\|_{E_{r_{0},1}}$. After that
to obtain \eqref{5.30.4} it suffices to use
again \eqref{9.28.1} with $\alpha=1$ and $\beta=1+\gamma p' $ to see that
$$
P_{1} \big(b^{ 1+\gamma p'}\big)\leq N\|b\|_{E_{ 1+\gamma p',1}}\big(  M\big(b^{ 1+\gamma p'}\big)\big)^{1-1/(1+\gamma p')}
\leq N\|b\| _{E_{r_{0},1}}b^{\gamma p'}.
$$

We now get rid of the assumption that 
$\hat M(|b|^{ r_{0} })\leq N_{0}|b|^{ r_{0 }}$ as in \cite{CF_90}. 

For $ r_{1} =( r_{0} +q)/2$ we have $|b|^{ r_{0} }\leq 
(\hat M(|b|^{ r_{1} }))^{ r_{0 }/  r_{1}} :=\tilde b^{ r_{0} }$ and since $ r_{0} / r_{1} <1$, $\tilde b^{  r_{0} }$ is an $A_{1}$-weight
with $N_{0}=N_{0}( r_{0} / r_{1} )$ (see, for instance, \cite{GR_85} p.~158).
Therefore, \eqref{5.30.6} holds with $\tilde b$
in place of $b$ and it only remains to show that  
$$
\|\tilde b\|_{E_{ r_{0} ,1}}\leq N\|b\|_{E_{q,1}},
$$
that is, for any $t,x,\rho$,
\begin{equation}
                             \label{1.19.1}
\int_{C_{\rho}(t,x)}\tilde b^{ r_{0} }\,dz
\leq  N\rho^{d+2- r_{0} }\|b\|_{E_{q,1}}^{ r_{0} }.
\end{equation}
Of course, we may assume that $t=0,x=0$.
Then by H\"older's inequality we see that
the left-hand side of \eqref{1.19.1} is less than
$$
N\rho^{(d+2)(q- r_{0} )/q}\Big(\int_{\bR^{d+1}}
(\hat M(|b|^{ r_{1} }))^{q/ r_{1} }I_{C_{\rho}}\,dz\Big)^{ r_{0 }/q},
$$
where the integral 
by a Fefferman-Stein Lemma 1, p.~111 of
\cite{FS_71} and the fact that
$q/ r_{1 }>1$ is dominated by
$$
N\int_{\bR^{d +1 }}|b|^{q}\hat MI_{C_{\rho}}\,dz
\leq N\rho^{d+2-q}\|b\|_{E_{q,1}}^{q},
$$
where we used Lemma \ref{lemma 1.22.1} b) for $\alpha=\beta=1$.
Hence, 
$$
\int_{C_{r} }\tilde b^{ r_{0} }\,dz
\leq N\rho^{(d+2)(q- r_{0} )/q+(d+2-q) r_{0} /q}\|b\|_{E_{q,1}}^{ r_{0} },
$$
which is \eqref{1.19.1}. 

An alternative way to get the result is
to follow the proof of Theorem 3 of \cite{CF_88}. We have
$$
\int_{\bR^{d+1}}
(\hat M(|b|^{ r_{1} }))^{q/ r_{1 }}I_{C_{\rho}}\,dz
\leq \int_{\bR^{d+1}}
(\hat M(|b|^{ r_{0} }))^{q/ r_{0} }(\hat MI_{C_{\rho}})^{\alpha}\,dz=:J,
$$
where $\alpha\in(0,1)$. An easy exercise
leads to the well-known result that
$(\hat MI_{C_{\rho}})^{\alpha}$ is an $A_{1}$-weight, and, hence, an  $A_{q/ r_{0} }$-weight.
By the Muckenhoupt theorem 
$$
J\leq N\int_{\bR^{d+1}}
 |b|^{q} (\hat MI_{C_{\rho}})^{\alpha}\,dz
$$
and it only remains to use 
Lemma \ref{lemma 1.22.1} b) again  with $ \beta=1$ and any appropriate $\alpha$.
The theorem is proved.

\mysection{Mixed-norm estimates}

                       \label{section 10.8.3}
For $q_{1},q_{2}\in[1,\infty]$ and measurable $f$ and $\Gamma\subset \bR^{d+1}$ introduce
$$
\| f\|_{L_{q_{1},q_{2}} } 
=\Big(\int_{\bR}\Big(\int_{\bR^{d}}|f(t,x)|^{q_{1}}
\,dx\Big)^{q_{2}/q_{1}}\,dt\Big)^{1/q_{2}},
$$
$$
\dashnorm f\|_{L_{q_{1},q_{2}}(\Gamma)}
 =\|  I_{\Gamma}\|_{L_{q_{1},q_{2}} }^{-1}\| fI_{\Gamma}\|_{L_{q_{1},q_{2}} }.
$$
Here the index of $L_{q_{1},q_{2}}$ which is the exponent of $\rho$ in the expression
$$
\| I_{C_{\rho}}\|_{L_{q_{1},q_{2}} }\quad
\text{\rm is}\quad \frac{d}{q_{1}}+\frac{2}{q_{2}}.
$$

If in addition  $0<\beta
\leq d/q_{1}+2/q_{2}$, set  
$$
\|f\|_{E_{q_{1},q_{2},\beta}(Q)}=
\sup_{\rho<\infty,(t,x)\in Q}
\rho^{\beta}\dashnorm I_{Q} f\|_{L_{q_{1},q_{2}}(C_{\rho}(t,x))}.
$$
We also introduce the spaces $L_{q_{1},q_{2}}(Q)$
and $E_{q_{1},q_{2},\beta}(Q)$ as the spaces of
functions whose respective norms are finite.
We abbreviate $L_{q_{1},q_{2}} =L_{q_{1},q_{2}}
( \bR^{d+1})$, $E_{q_{1},q_{2},\beta} =E_{q_{1},q_{2},\beta}
( \bR^{d+1})$.

The following is certainly well known.

\begin{lemma}
                         \label{lemma 10.6.1}
Let $f$ be a nonnegative function on $\bR^{d+1}$, $p,q\in(1,\infty)$. Then for
any $w_{x}(x),w_{t}(t)$ which are $A_{p}$
Muckenhoupt weights on $\bR^{d}$ and $\bR$,
respectively, we have
\begin{equation}
                             \label{10.6.3}
\int_{\bR^{d+1}}|\hat Mf|^{p}w_{x}w_{t}
\,dxdt\leq N\int_{\bR^{d+1}}| f|^{p}w_{x}w_{t}
\,dxdt,
\end{equation}
where $N$ depends only on $d,p$, and the $A_{p}$-constants of $w_{x},w_{t}$. Furthermore,
\begin{equation}
                             \label{10.6.4}
\int_{-\infty}^{\infty}\Big(\int_{\bR^{d}}
|\hat Mf|^{p}\,dx\Big)^{q/p}\,dt
\leq N \int_{-\infty}^{\infty}\Big(\int_{\bR^{d}}
| f|^{p}\,dx\Big)^{q/p}\,dt,
\end{equation}
where $N$  depends only on $d,p,q$.

\end{lemma}

Proof. Estimate \eqref{10.6.3} follows
by  application of the Muckenhoupt
theorem to $w_{x}w_{t}$, which is an
$A_{p}$-weight on $\bR^{d+1}$. Then observe that in the particular case that $w_{x}\equiv 1$, \eqref{10.6.3} means
that
$$
\int_{-\infty}^{\infty}\Big[\Big(\int_{\bR^{d}}
|\hat Mf|^{p}\,dx\Big)^{1/p}\Big]^{p}w_{t}\,dt
\leq N \int_{-\infty}^{\infty}\Big[\Big(\int_{\bR^{d}}
| f|^{p}\,dx\Big)^{1/p}\Big]^{p}w_{t}\,dt
$$
for any $A_{p}$-weight $w_{t}$, which implies \eqref{10.6.4} by the Rubio de Francia
extrapolation theorem. The lemma is proved.

This lemma, \eqref{9.28.1}, and H\"older's inequality immediately yield the following.

\begin{lemma}
                           \label{lemma 10.7.1}
For any $\alpha\in(0,\beta ),\beta\in(0,d+2]$, 
 $p\in[1,\infty]$,
$q_{1},q_{2}\in(1,\infty]$, there exists a constant $N$ such that for any $f\geq0$    and measurable $\Gamma$ we have
\begin{equation}
                                \label{10.7.1}
 \|P_{\alpha}f\|_{L_{r_{1},r_{2}}(\Gamma)} \leq    
N\|M_{\beta}f \|_{L_{p}(\Gamma)}^{\alpha/\beta}
\|f\|_{L_{q_{1},q_{2}}}^{1-\alpha/\beta},
\end{equation}
provided that
$$
\frac{1}{r_{i}}=\frac{\alpha}{\beta}\cdot\frac{1}{p}+
\Big(1-\frac{\alpha}{\beta}\Big)\frac{1}{q_{i}},\quad i=1,2.
$$
\end{lemma}

Similarly to Corollary \ref{corollary 9.29.1}
we have
 \begin{corollary}
                    \label{corollary 10.7.1}
Let    
$q_{1},q_{2}\in(1,\infty]$,
$$
\beta:=\frac{d}{q_{1}}+\frac{2}{q_{2}}>0,
$$
$\alpha\in(0,\beta )$. Then for any $f\geq0$
we have
$$
\|P_{\alpha}f\|_{L_{r_{1},r_{2}}}\leq N\|f\|_{L_{q_{1},q_{2}}}
$$
as long as $q_{i}\beta=r_{i}(\beta-\alpha)$,
$i=1,2$.

In particular, (almost follows from
Theorem 10.2 of \cite{BIN_75})
if    $\beta>1$, and    $u\in C^{\infty}_{0} $, then
\begin{equation}
                            \label{10.8.4}
\|Du\|_{L_{r_{1},r_{2}}}\leq N\|\partial_{t}u+\Delta u\|_{L_{q_{1},q_{2}}}
\end{equation}
as long as
$q_{i}\beta=r_{i}(\beta-1)$,
$i=1,2$.
\end{corollary}

\begin{corollary}
                     \label{corollary 10.8.2}
Under the assumptions of Corollary
\ref{corollary 10.7.1}, if $\beta>1$, there is a constant $N$ such that, for any $b=(b^{i})
\in L_{\beta q_{1} ,\beta q_{2}}$ and $u\in C^{\infty}_{0} $,
\begin{equation}
                         \label{2.29.2}
\|b^{i}D_{i}u\|_{L_{q_{1},q_{2}}}\leq N
\|b  \|_{L_{\beta q_{1},\beta q_{2}}}
\|\partial_{t}u+\Delta u\|_{L_{q_{1},q_{2}}}.
\end{equation}

\end{corollary}

Indeed, by H\"older's inequality
$$
\|b^{i}D_{i}u\|_{L_{q_{1},q_{2}}}
\leq \|b  \|_{L_{\beta q_{1},\beta q_{2}}}
\|Du\|_{L_{r_{1},r_{2}}}.
$$
\begin{remark}
                      \label{remark 10.8.30}
It is instructive to compare this result
with Remark \ref{remark 10.8.5}. Now
we can treat $b\in L_{s_{1},s_{2}}$ with
$s_{i}\in(1,\infty]$ satisfying $d/s_{1}
+2/s_{2}=1$.
\end{remark}

Since $E_{q_{1},q_{2},\beta}=L_{q_{1},q_{2}}$
if $\beta=d/q_{1}+2/q_{2}$,
the following is a generalization of
Corollary \ref{corollary 10.7.1}.

\begin{theorem}
                       \label{theorem 10.7.1}
Let    
$q_{1},q_{2}\in(1,\infty]$,
$$
\frac{d}{q_{1}}+\frac{2}{q_{2}}\geq \beta>0,
$$
$\alpha\in(0,\beta )$. Then 
there is a constant $N$ such that for any $f\geq0$
we have
\begin{equation}
                          \label{10.7.4}
 \|P_{\alpha}f\|_{E_{r_{1},r_{2},\beta-\alpha}}
\leq N \|f\|_{E_{q_{1},q_{2},\beta}},
\end{equation}
where $r_{i}(\beta-\alpha)=q_{i}\beta$,
$i=1,2$.
\end{theorem}

Proof. It suffices to prove that for any $\rho>0$
$$
\rho^{\beta-\alpha}\Big(\dashint_{0}^{\,\,\,\rho^{2}}\Big(\dashint_{B_{\rho}}|P_{\alpha}f|^{r_{1}}\,dy\Big)^{r_{2}/r_{1}}\,ds\Big)^{1/r_{2}}
\leq N\|f\|_{E_{q_{1},q_{2},\beta}},
$$
that is  
\begin{equation}
                          \label{10.7.6}
\rho^{\beta-\alpha-(d/r_{1}+2/r_{2})}\Big(\int_{0}^{ \rho^{2}}\Big(\int_{B_{\rho}}|P_{\alpha}f|^{r_{1}}\,dy\Big)^{r_{2}/r_{1}}\,ds\Big)^{1/r_{2}}
\leq N\|f\|_{E_{q_{1},q_{2},\beta}}.
\end{equation}

Observe that by H\"older's inequality
$M_{\beta}f\leq N\|f\|_{E_{q_{1},q_{2},\beta}}$ and by definition
$$
\|I_{C_{2\rho}}f\|_{L_{q_{1},q_{2}}}=
N\rho^{d/q_{1}+2/q_{2}}
\Big(\dashint_{0}^{\,\,\,4\rho^{2}}\Big(\dashint_{B_{2\rho}}| f|^{q_{1}}\,dy\Big)^{q_{2}/q_{1}}\,ds\Big)^{1/q_{2}}
$$
$$
\leq N\rho^{d/q_{1}+2/q_{2}-\beta} \|f\|_{E_{q_{1},q_{2},\beta}}
=N\rho^{(d/r_{1}+2/r_{2})\beta/(\beta-\alpha)-\beta} \|f\|_{E_{q_{1},q_{2},\beta}}.
$$
 
It follows from  Lemma \ref{lemma 10.7.1} 
with $p=\infty$ that \eqref{10.7.6} holds with
$I_{C_{2\rho}}f$ in place of $f$ on the left.
 
Furthermore,  by Corollary \ref{corollary 1.17.1} we have 
 $|P_{\alpha}(I_{C_{2\rho}^{c}}f)|
\leq N\rho^{\alpha-\beta}M_{\beta}f$ in $C_{\rho}$.
Therefore,
$$
\rho^{\beta-\alpha}\Big(\dashint_{0}^{\,\,\, \rho^{2}}\Big(\dashint_{B_{\rho}}|P_{\alpha}(I_{C_{2\rho}^{c}}f)|^{r_{1}}\,dy\Big)^{r_{2}/r_{1}}\,ds\Big)^{1/r_{2}}\leq N\sup M_{\beta}f
\leq N\|f\|_{E_{q_{1},q_{2},\beta}}.
$$
By combining these results we come to 
\eqref{9.29.4} and the theorem is proved.

\begin{corollary}
                      \label{corollary 10.8.1}
Under the assumptions of Theorem \ref{theorem 10.7.1}, if $\beta>1$,  for any $u\in C^{\infty}_{0}$, we have
$$
\|Du\|_{E_{r_{1},r_{2},\beta-1}}
\leq N\|\partial_{t}u+\Delta u\|_{E_{q_{1},q_{2},\beta }},
$$
where $r_{i}(\beta-1)=q_{i}\beta$,
$i=1,2$. This coincides with \eqref{10.8.4} if
$\beta$ is equal to the index of $L_{q_{1},q_{2}}$.
\end{corollary}

\begin{remark}
                      \label{remark 10.8.2}
Corollary \ref{corollary 10.8.1} opens up the
possibility to treat the terms like
$b^{i}D_{i}u$ as perturbation terms
in operators like $\partial_{t}u+\Delta u
+b^{i}D_{i}u$ with even  lower
summability properties of $b=(b^{i})$
than in Remark \ref{remark 10.8.30}.
To show this observe that  for
$q_{1},q_{2},\beta$ as in Theorem \ref{theorem 10.7.1} with $\beta>1$ and
$s_{i}=\beta q_{i}\in(1,\infty]$, $i=1,2$, we have
$$
\rho^{\beta}\dashnorm I_{C_{\rho}}b^{i}D_{i}u\|_{L_{q_{1},q_{2}}}
\leq \rho\dashnorm b I_{C_{\rho}}\|_{L_{s_{1},s_{2}}}\cdot \rho^{\beta-1}
\dashnorm I_{C_{\rho}} Du \|_{L_{r_{1},r_{2}}}
$$
implying that
\begin{equation}
                         \label{2.29.1}
\| b^{i}D_{i}u\|_{E_{q_{1},q_{2},\beta}}
\leq \|b  \|_{E_{s_{1},s_{2},1}} \|   Du  \|_{E_{r_{1},r_{2},\beta-1}}\leq N\|b  \|_{E_{s_{1},s_{2},1}}
\|\partial_{t}u+\Delta u\|_{E_{q_{1},q_{2},\beta }},
\end{equation}
where  
$d/s_{1}+2/s_{2}\geq 1$.

However, note that we also need
$$
\rho\dashnorm b I_{C_{\rho}(t,x)}\|_{L_{s_{1},s_{2}}}
$$
to be bounded as a function of $\rho, t,x$.
If we ask ourselves what $\tau>0$ should be
to guarantee this boundedness if
$b\in L_{\tau s_{1},\tau s_{2}}$, if
$d/s_{1}+2/s_{2}>1$,
the slightly disappointing answer is that
$\tau=d/s_{1}+2/s_{2}$, so that
$d/(\tau s_{1})+2/(\tau s_{2})=1$. Still
functions in $E_{s_{1},s_{2},1}$ may have higher singularities than those in
$L_{\tau s_{1},\tau s_{2}}$.

 Another advantage of 
\eqref{2.29.1} in comparison with \eqref{2.29.2} is seen when $b$ depends only
on $t$ or $|b(t,x)|\leq \hat b(t)$. In that case \eqref{2.29.1}
becomes 
$$
\| b^{i}D_{i}u\|_{E_{q_{1},q_{2},\beta}}
\leq N\|\hat b  \|_{E_{ \beta q_{2},1/2}(\bR)}
\|\partial_{t}u+\Delta u\|_{E_{q_{1},q_{2},\beta }},
$$
and if $\beta q_{2}=2$, then
$$
\|\hat b  \|_{E_{ \beta q_{2},1/2}(\bR)}=\|\hat b\|_{L_{2}(\bR)}.
$$
Thus for any $q_{1}\in(1,\infty]$ and  $q_{2}\in(1,2)$
$$
\| b^{i}D_{i}u\|_{E_{q_{1},q_{2},2/q_{2}}}
\leq N \|\hat b\|_{L_{2}(\bR)}
\|\partial_{t}u+\Delta u\|_{E_{q_{1},q_{2},2/q_{2}}}.
$$
In case  $q_{1}\in(1,d)$, $q_{2}\in(1,\infty]$, $1<\beta\leq d/q_{1}$, and $\|b(\cdot,t)\|_{E_{\beta q_{1},1}(\bR^{d})}\leq \hat b<\infty$ for any $t$,
we also have
$$
\| b^{i}D_{i}u\|_{E_{q_{1},q_{2},\beta}}
\leq N\hat b
\|\partial_{t}u+\Delta u\|_{E_{q_{1},q_{2},\beta }}.
$$
An application of the last inequality
in case $u,b$ are independent of $t$,
  $\beta=d/q_{1}$, $q_{1}\in(1,d)$, and $q_{2}=\infty$, yields the well-known estimate
$$
 \| b^{i}D_{i}u \|_{L_{q_{1}}(\bR^{d})}\leq N \| b \|_{L_{d}(\bR^{d})}
 \|  \Delta u \|_{L_{q_{1}}(\bR^{d})}.
$$
 
\end{remark}

To extend the embedding and interpolation results
to Morrey spaces with mixed norms
we need the following  result
very useful also in other circumstances.

\begin{lemma}[Poincar\'e's inequality]
                     \label{lemma 10.10.1}
Let $1\leq r_{1},r_{2}<\infty$,
$u\in C^{\infty}_{0}$, $\rho\in(0,\infty)$.
Then
\begin{equation}
                         \label{10.10.1}
\dashnorm Du-(Du)_{C_{\rho}}\|_{L_{r_{1},r_{2}}(C_{\rho})}^{r_{2}}\leq
N(d,r_{1},r_{2})\rho^{r_{2}}\dashnorm \,|\partial_{t}u|+|D^{2}u| \|_{L_{r_{1},r_{2}}(C_{\rho})}^{r_{2}}.
\end{equation}

\end{lemma}

Proof. We follow the usual way (see, for instance,
Lemma 4.2.2 of \cite{Kr_08}). First, due
to self-similar transformations, we may take
$\rho=1$. In that case,  for
a $\zeta\in C^{\infty}_{0}(B_{1})$ with unit integral, introduce
$$
v(t)=\int_{B_{1}}\zeta(y)Du(t,y)\,dy.
$$
Then by the usual Poincar\'e inequality
$$
\int_{B_{1}}|Du(t,x)-v(t)|^{r_{1}}\,dx=
\int_{B_{1}}\big|\int_{B_{1}}[Du(t,x)-Du(t,y)]\zeta(y)\,dy\big|^{r_{1}}\,dx
 $$
 \begin{equation}                
                              \label{7.8.3}
\leq N\int_{B_{1}}\int_{B_{1}}
|Du(t,x)-Du(t,y)|^{r_{1}}\,dxdy
\leq N\int_{B_{1}}|D^{2}u(t,x)|^{r_{1}}\,dx.
\end{equation}

Next, observe that for any constant vector $v$
the left-hand side of
\eqref{10.10.1} is less than a constant times (recall that $\rho=1$)
 $$
\int_{0}^{1}\Big(\int_{B_{1}}|Du(t,x)-v|^{r_{1}}\,dx\Big)^{r_{2}/r_{1}}dt
$$
$$
\leq N
\int_{0}^{1}\Big(\int_{B_{1}}|Du(t,x)-v(t)|^{r_{1}}\,dx\Big)^{r_{2}/r_{1}}dt+N
\int_{0}^{1}|v(t)-v|^{r_{2}}\,dt.
 $$
By \eqref{7.8.3} the first term on the right is less
than the right-hand side of \eqref{10.10.1}. 
To estimate the second term,
take 
 $$
v=\int_{0}^{1}v(t)\,dt.
 $$
Then by Poincar\'e's inequality
 $$
\int_{0}^{1}|v(t)-v|^{r_{2}}\,dt\leq N
\int_{0}^{1}\big|\int_{B_{1}}\zeta 
\partial_{t} Du\,dx\big|^{r_{2}}\,dt
=N
\int_{0}^{1}\big|\int_{B_{1}}(D\zeta) 
\partial_{t}u\,dx\big|^{r_{2}}\,dt
$$
and to finish the proof it only remains
to use  
 H\"older's inequality.
 The lemma is proved.

The usual Poincar\'e inequality was used in the proof
of Lemma \ref{lemma 7.30.1}.
Also observe that mixed-norms estimates
like \eqref{7.31.1} are available in
\cite{BIN_75} (see Theorem 9.5 there).
Therefore, by using Lemma \ref{lemma 10.10.1}
and  following very closely the proofs of 
Lemmas \ref{lemma 10.4.1}, \ref{lemma 7.30.1}, and
Theorems \ref{theorem 10.2.2}  
we arrive at the following
results about interpolation 
and embedding for
  Morrey spaces with mixed norms.

\begin{lemma}
                      \label{lemma 10.10.2}
Let   $q_{1},q_{2}\in(1,\infty)$, $0<\beta\leq d/q_{1}+2/q_{2}$. Then there is a constant
$N $ such that, for any 
$R\in(0,\infty)$, 
   $\varepsilon
\in(0,1] $,  
and $u\in C^{\infty}_{0}$,
 \begin{equation}
                              \label{10.10.4}
\|Du\|_{E_{q_{1},q_{2},\beta}(C_{R})}\leq
N\varepsilon R \||\partial_{t}u|+|D^{2}u|\,\|_{E_{q_{1},q_{2},\beta}(C_{R})}
+N\varepsilon^{-1}R^{-1}
\| u\|_{E_{q_{1},q_{2},\beta}(C_{R})}.
\end{equation}
\end{lemma}

\begin{theorem}
                        \label{theorem 10.10.2}
Let   $q_{1},q_{2}\in(1,\infty)$, $1<\beta\leq d/q_{1}+2/q_{2}$ and let $r_{i}(\beta-1)=q_{i}\beta$, $i=1,2$.  
 Then  there is a constant $N$ such that for any $R\in(0,\infty]$, $u\in C^{\infty}_{0}$
we have
\begin{equation}
                                  \label{10.10.2}
 \|Du\|_{E_{r_{1},r_{2},\beta-1}(C_{R})}\leq N\||\partial_{t}u|+|D^{2} u|\,\|_{E_{q_{1},q_{2},\beta}(C_{R})}
+NR^{-2}\|u\|_{E_{q_{1},q_{2},\beta}(C_{R})}.
\end{equation}

\end{theorem}
 
\begin{remark}
By taking $u$ depending only on $x$ we recover from Lemma \ref{lemma 10.10.2} and Theorem \ref{theorem 10.10.2} their ``elliptic''
counterpart stated as Lemmas 4.4 
and 4.7 in \cite{Kr_303}, respectively.
\end{remark}

\begin{remark}
  Theorem \ref{theorem 10.10.2} is the most general results of the paper containing as particular cases our previous results on  embeddings. Thus, Corollary \ref{corollary 10.8.1} (in an obvious rougher form) follows from Theorem \ref{theorem 10.10.2} when $R=\infty$ and contains embedding results
for Lebesgue spaces with mixed norms as $\beta
=d/q_{1}+2/q_{2}$ and for $L_{q}$-spaces as $q=q_{1}=q_{2}$.
\end{remark}

\begin{remark}
                       \label{remark 10.10.4}
We stated our results only for $u\in C^{\infty}_{0}$ just for convenience. Let us show why, for instance, 
Theorem \ref{theorem 10.10.2} is valid as long as $\partial_{t}u,Du,D^{2}u\in E_{q_{1},q_{2},\beta}(C_{R})$. For that, it suffices to prove that for any $R'<R$, $\rho>0$, $(t,x)
\in C_{R'}$ the quantity
$$
I:=\rho^{\beta}\dashnorm
 I_{C_{R'}} Du\|_{L_{r_{1},r_{2}}(C_{\rho}(t,x))}
$$
is less than the right-hand side of \eqref{10.10.2} with $(R')^{-2}$
in place of $R^{-2}$. For $\varepsilon>0$ define $u^{(\varepsilon)}=
(I_{C_{R}}u)*\zeta_{\varepsilon}$, where $\zeta_{\varepsilon}(x)=\varepsilon^{-d-1}
\zeta(t/\varepsilon,x/\varepsilon)$, nonnegative
$\zeta\in C^{\infty}_{0}$ has integral one and $\zeta(t,x)=0$
for $t\geq0$. Also  
introduce $I^{\varepsilon}$ by replacing $u$
in the definition of $I$ with $u^{(\varepsilon)}$. Of course,
$I^{\varepsilon}\to I$ as $\varepsilon
\downarrow 0$ and by Theorem \ref{theorem 10.10.2}
$$
I^{\varepsilon}\leq 
N\||\partial_{t}u^{(\varepsilon)}|+|D^{2} u^{(\varepsilon)}|\,\|_{E_{q_{1},q_{2},\beta}(C_{R'})}
+N(R')^{-2}\|u^{(\varepsilon)}\|_{E_{q_{1},q_{2},\beta}(C_{R'})}=:J^{\varepsilon}.
$$

Observe that if $\varepsilon$ is small enough
and $(s,y)\in C_{R'}$, then $\partial_{t}
u^{(\varepsilon)}(s,y)=(I_{C_{R}}\partial_{t}u)
*\zeta_{\varepsilon}(s,y)$. Similar formulas
are valid for $D^{2}u^{(\varepsilon)}$ and by
Minkowski's inequality (the norm of a sum
is less then the sum of norms) we have
$$
J^{\varepsilon}\leq  \int_{\bR^{d+1}}\zeta(s,y)\Big(N\|I_{C_{R}}(|\partial_{t}u |+|D^{2} u |)(\cdot-\varepsilon(s,y)\,\|_{E_{q_{1},q_{2},\beta}(C_{R'})}
$$
$$
+N(R')^{-2}\|I_{C_{R}}u (\cdot-\varepsilon(s,y) \|_{E_{q_{1},q_{2},\beta}(C_{R'})}\Big)\,dyds
$$
$$
= \int_{\bR^{d+1}}\zeta(s,y)\Big(N\|I_{C_{R}}(|\partial_{t}u |+|D^{2} u |) \,\|_{E_{q_{1},q_{2},\beta}(C_{R'}-\varepsilon(s,y))}
$$
$$
+N(R')^{-2}\|I_{C_{R}}u   \|_{E_{q_{1},q_{2},\beta}(C_{R'}-\varepsilon(s,y))}\Big)\,dyds.
$$
Since in the last integral $C_{R'}-\varepsilon(s,y)\subset C_{R}$ if $\varepsilon$
is small enough, it follows that
for small $\varepsilon$
$$
J^{\varepsilon}\leq 
N\|I_{C_{R}}(|\partial_{t}u |+|D^{2} u |) \,\|_{E_{q_{1},q_{2},\beta}(C_{R })}
+N(R')^{-2}\|I_{C_{R}}u   \|_{E_{q_{1},q_{2},\beta}(C_{R } )}
$$
which yields the desired result.
\end{remark}

{\bf Acknowledgments}. The author brings his gratitude to the referees for their useful comments.

\end{document}